\documentclass{amsart}
\usepackage{amscd} 
\usepackage{mathrsfs} 
\usepackage{amssymb}

\theoremstyle{plain} 
\newtheorem{thm}{Theorem}[section]
\newtheorem{lem}[thm]{Lemma} 
\newtheorem{prop}[thm]{Proposition}
\newtheorem{cor}[thm]{Corollary}

\theoremstyle{definition} 
\newtheorem{defn}[thm]{Definition}
\newtheorem{exmp}[thm]{Example}

\theoremstyle{remark}
\newtheorem{rem}[thm]{Remark}
\newtheorem*{ackn}{Acknowledgment}

\newcommand{\gCi}{C^\infty}
\newcommand{\tans}{{\mathcal{T}}}
\newcommand{\tang}{{\mathcal{F}}}
\newcommand{\trans}{{\mathcal{Q}}}
\newcommand{\cl}[1]{[{#1}]}
\newcommand{\field}[1]{\mathbb{#1}}
\newcommand{\R}{\field{R}}
\newcommand{\C}{\field{C}}
\newcommand{\Z}{\field{Z}}
\newcommand{\gSec}{\varGamma}
\newcommand{\pt}{\odot}

\DeclareMathOperator{\ctens}{\hat{\otimes}}
\DeclareMathOperator{\id}{id}
\DeclareMathOperator{\im}{im}
\DeclareMathOperator{\Supp}{Supp}
\DeclareMathOperator{\sgn}{sgn}
\DeclareMathOperator{\pr}{pr}
\DeclareMathOperator{\Coin}{Coin}
\DeclareMathOperator{\vol}{vol}
\DeclareMathOperator{\Tr}{Tr}

\begin{document}

\title[Coincidence Formula]{A Coincidence Formula for Foliated
  Manifolds}

\author{Bernd M\"umken}

\address{Mathematisches Institut, West\-f\"alische
  Wilhelms-Universit\"at M\"unster, Einstein\-str.~62,
  48149~M\"unster, Germany.}

\email{mumken@math.uni-muenster.de}


\begin{abstract}
  The main result of the present paper is a coincidence formula for
  foliated manifolds. To prove this we establish K\"unneth formula,
  Poincar\'e duality and intersection product in the context of
  tangential de~Rham cohomology and homology of tangential currents.
  
  We apply the formula to get a dynamical Lefschetz formula for
  foliated flows.
\end{abstract}

\maketitle{}

\section{Introduction}
\label{sec:intro}

\subsection{{}}
\label{sec:introa}

Consider the category of foliated manifolds with foliated maps. A
foliated manifold~$M$ is a smooth, i.e.~$\gCi{}$, manifold without
boundary together with a decomposition into connected subsets, called
leaves, which is locally trivial. Here local triviality means that
each point of~$M$ has an open neighborhood~$U$ which is diffeomorphic
to a Euclidean vector space~$V$ such that the induced partition on~$U$
corresponds to the partition of~$V$ into the cosets of a vector
subspace~$W$ of~$V$. The quotient~$V/W$ is called a local transversal
manifold. A foliated map is a smooth map such that leaves are mapped
to leaves.

Denote by~$\tans{}$ the covariant tangent functor. Attach to a foliated
manifold~$M$ the tangential subbundle~$\tang{}M$ of~$\tans{}M$
consisting of vectors which are tangent to the leaves. Since the
derivative of a foliated map induces a map of tangential subbundles we
get a functor~$\tang{}$. Passing to quotients brings the functor
$\trans{}=\tans{}/\tang{}$. For a foliated manifold~$M$ the vector
bundle~$\trans{}M$ is the transversal bundle of~$M$.

Let $M$ be a foliated manifold. A transversal Riemannian metric on~$M$
is a Riemannian metric in~$\trans{}M$ which is locally the pull back
of a metric on the local transversal manifold. Next identify
$\tans{}M=\tang{}M\oplus{}\trans{}M$ by a splitting. Then a bundle-like
metric on~$M$ is a Riemannian metric in~$\tans{}M$ which is the direct
sum of a Riemannian metric in~$\tang{}M$ and a transversal Riemannian
metric on~$M$.

Mimic the proceeding in the classical de~Rham theory with the whole
tangent bundle replaced by the tangential subbundle to attain to the
contravariant functors $\Omega^{*}_{\tang{}}$
and~$\Omega^{*}_{\tang{},c}$, the tangential differential forms (with
compact support). These are functors to nuclear graded differential
$\R$-vector spaces. The corresponding cohomology functors
$H^{*}_\tang{}$ and~$H^{*}_{\tang{},c}$ are the tangential cohomology
and the tangential cohomology with compact support.

Equip the topological dual~$E'$ of a topological vector space~$E$ with
the weak topology. Applying the functor continuous linear forms to
tangential differential forms define tangential currents with compact
support~$\Omega^{\tang{}}_{c,*}$. Similarly, there are the tangential
currents~$\Omega^{\tang{}}_{*}$ dual to the tangential differential forms
with compact support. In this way, find the covariant functors
$H^\tang_{c,*}$ and~$H^\tang_{*}$, the tangential homology (with
compact support).

The differentials just introduced are not topological homomorphisms in
general. Consequently, it is reasonable to consider the maximal
Hausdorff quotients of tangential (co)homology (with compact support).
This reduction preserves the functorial properties. The resulting
functors are denoted by $\bar{H}^{*}_\tang$,~$\bar{H}^{*}_{\tang,c}$,
$\bar{H}^\tang_{c,*}$ and~$\bar{H}^\tang_{*}$ and are called reduced
tangential (co)homology (with compact support). In particular, get a
canonical topological isomorphism
\begin{displaymath}
  \bar{H}^{\tang}_{*}(M)\cong\bar{H}^{*}_{\tang,c}(M)'
\end{displaymath} 
for each foliated manifold~$M$.

Let $M$ be a foliated manifold. Then the zeroth reduced tangential
cohomology of~$M$ (with compact support) consists of those smooth
functions (with compact support) which are constant along the leaves.
Now let $M$ carry a transversal Riemannian metric and assume that the
tangential subbundle of~$M$ is oriented. Let the leaves of~$M$ have
dimension~$p$. Then a non trivial reduced homology class is defined by
integration of tangential $p$-forms against the transversal
volume~$\vol_\trans$ induced by the metric. This is denoted by
$\int_M-\:{}\vol_\trans$.

Now we present the main results of this paper in the setting of the
motivating application. In this context our coincidence formula is
called a dynamical Lefschetz formula.

\subsection{{}} 
\label{sec:introb}

Let the real line~$\R$ be foliated by points and equipped with the
standard metric and orientation. Let $M$ be a compact connected
foliated manifold with a bundle-like metric and oriented tangential
subbundle. Let $p$ be the dimension of the leaves of~$M$.

Denote by $-\ctens{}-$ the complete tensor product of nuclear vector
spaces. Then we have a K\"unneth theorem,
cf.~Section~\ref{sec:kuenneth}.

\begin{thm} 
  \label{thm:kuenneth} 
  The canonical maps
  \begin{displaymath}
    \gCi(\R_{>0})\ctens\bar{H}_\tang^{*}(M)\to\bar{H}^{*}_{\tang}(\R_{>0}\times{}M) 
  \end{displaymath}
  and
  \begin{displaymath}
    \bigoplus_\kappa\gCi_c(\R_{>0})\ctens\bigl(\bar{H}_\tang^{p-\kappa}(M)\ctens\bar{H}_\tang^{\kappa}(M)\bigr)\to\bar{H}^{p}_{\tang,c}(\R_{>0}\times{}M\times{}M)
  \end{displaymath} 
  are topological isomorphisms.
\end{thm} 

The first statement is similar to a classical result concerning smooth
functions with values in a nuclear Fr\'echet space. The proof of the second
involves a tangential Hodge decomposition established by \'Alvarez
L\'opez and Kordyukov~\cite[Corollary~$1.3$]{MR1815391}.

Special cases of Poincar\'e duality take the following shape in this
context, cf.~Theorem~\ref{thm:regularb}.

\begin{thm} 
  \label{thm:poincareduality} 
  The maps
  \begin{align*}
    \bar{H}^{p-*}_\tang(M)&\to\bar{H}^{*}_{\tang}(M)',\\
    \omega&\mapsto\int_{M}\omega\wedge-\:\vol_\trans,
  \end{align*}
  and the similar map
  \begin{displaymath}
    \bar{H}^{p}_\tang(\R_{>0}\times{}M\times{}M)\to\bar{H}^{p}_{\tang,c}(\R_{>0}\times{}M\times{}M)'  \end{displaymath}
  are injective and have dense image.
\end{thm}

For finite dimensional $\bar{H}^{*}_\tang(M)$ the first map is an
isomorphism. Here the approximation is done by generalizing the
convolution of a distribution with a regularizing sequence of smooth
functions. More precisely, for the second result we define a sequence
of maps
\begin{displaymath}
  R_\nu:\bar{H}^{p}_{\tang,c}(\R_{>0}\times{}M\times{}M)'\to\bar{H}^{p}_{\tang}(\R_{>0}\times{}M\times{}M),
\end{displaymath} 
such that
$\int_{\R_{>0}\times{}M\times{}M}R_\nu{}S\wedge-\:\vol_\trans$
converges to~$S$ for each continuous linear form
$S\in\bar{H}^{p}_{\tang,c}(\R_{>0}\times{}M\times{}M)'$.

\subsection{{}} 
\label{sec:introc}

Let
\begin{displaymath}
  f:\R_{>0}\times{}M\to{}M
\end{displaymath} 
be a foliated map. Then the graph of~$f$ is a foliated submanifold of
$\R_{>0}\times{}M\times{}M$ which is isomorphic to~$\R_{>0}\times{}M$. Moreover, the K\"unneth theorem allows to define a reduced
tangential homology class $\cl{f}\in\bar{H}^\tang_p(\R_{>0}\times{}M\times{}M)$ by the requirement that
\begin{displaymath}
   \langle \cl{f}, \omega \otimes \tau \rangle = \int_{\R_{> 0} \times M} \omega \wedge f^* \tau \: \vol_\trans
\end{displaymath}
for $\omega \in \bar{H}^{t, p - \kappa}_c (\R_{> 0} \times M)$, $\tau
\in \bar{H}^{t, \kappa} (M)$ and $0 \le \kappa \le p$.

In particular, let now
\begin{align*}
  \phi :  \R \times M & \to M , \\
  (t, a) & \mapsto \phi^t (a)
\end{align*}
be a global foliated flow on~$M$ with generated vector field~$X$ and
let
\begin{displaymath}
  \pi : \R_{> 0} \times M \to M 
\end{displaymath} 
be the projection.

Suppose that the graphs $\Gamma$ of~$\phi$ and $\tilde{\Delta}$
of~$\pi$ intersect foliatedly transversally, i.e.
\begin{align*}
  \tans_{(t, a, a)} \Gamma + \tans_{(t, a, a)} \tilde{\Delta} & = \tans_{(t, a, a)} (\R_{> 0} \times M \times M), \\
  \tang_{(t, a, a)} \Gamma + \tang_{(t, a, a)} \tilde{\Delta} & =
  \tang_{(t, a, a)} (\R_{> 0} \times M \times M)
\end{align*} 
for all $(t, a) \in \R_{> 0} \times M$ with $\phi^t (a) = a$. Then the
intersection $\Gamma \cap \tilde{\Delta}$ is either empty or a one
dimensional submanifold of $\R_{> 0} \times M \times M$ which is
transversal to the leaves. It consists of components
\begin{displaymath}
  \R_{> 0} \times \{a\} \times \{a\} ,   
\end{displaymath}
where $a$ is a fixed point of~$\phi$, and
\begin{displaymath}
  \{\nu l (\gamma)\} \times \bigl\{ \bigl( \phi^t (b_\gamma) , \phi^t (b_\gamma) \bigr) \mid t \in \R \bigr\} ,
\end{displaymath}
$\nu = 1, 2 \dots$, where $\gamma$ is a periodic orbit with least
positive period~$l (\gamma)$ and $b_\gamma \in \gamma$ is arbitrary.
For each~$\gamma$ and~$\nu$ the derivative of~$\phi^{\nu l (\gamma)}$
induces an endomorphism
\begin{displaymath}
  \overline{\trans}_{b_\gamma} \phi^{\nu l (\gamma)} : \trans_{b_\gamma} M / \R \cdot \overline{X}_{b_\gamma} \to \trans_{b_\gamma} M / \R \cdot \overline{X}_{b_\gamma}.
\end{displaymath}
If we let $\pr : \R_{> 0} \times M \times M \to \R_{> 0}$ be the
projection then from Example~\ref{exmp:flowrestriction} we get

\begin{thm} 
  \label{thm:intersectionproduct} 
  The intersection product
  \begin{displaymath}
    \cl{\phi} \bullet \cl{\pi} := \lim_{\nu \to \infty} \langle \cl{\phi} , R_\nu \cl{\pi} \wedge - \rangle
  \end{displaymath}
  is defined as a continuous linear form on $\bar{H}^0_{\tang, c}(\R_{> 0} \times M \times M)$. Its push forward $\pr_* (\cl{\phi}
  \bullet \cl{\pi})$ is a distribution on~$\R_{>0}$ which is a sum of
  local contributions coming from the fixed points and the periodic
  orbits of~$\phi$. Explicitly,
  \begin{multline}
    \label{eq:intersectionproduct} 
    \pr_* (\cl{\phi} \bullet \cl{\pi}) = \sum_a \sgn \det (\id - \tang_{a} \phi^{t_a}) \int_{\R_{> 0}} \frac{1}{\bigl\vert \det \bigl( \id - \trans_a \phi^t \bigr) \bigr\vert} \cdot - \: d t \\
    + \sum_\gamma l (\gamma) \sum_{\nu = 1}^\infty \sgn \det (\id -
    \tang_{b_\gamma} \phi^{\nu l (\gamma)}) \frac{1}{\bigl\vert \det
      \bigl( \id - \overline{\trans}_{b_\gamma} \phi^{\nu l (\gamma)}
      \bigr) \bigr\vert} \cdot \delta_{\nu l (\gamma)} ,
  \end{multline} 
  where $ t_a > 0 $ is arbitrary and $\delta_{t_0}$~is the Dirac
  distribution in~$t_0$.
\end{thm}

The main step in the proof of this theorem is a direct computation of
the approximately defined intersection product.

\subsection{{}} 
\label{sec:introd}

Again let $\phi$ be a global foliated flow on~$M$ with generated
vector field~$X$ and let $\pi$ and~$\pr$ be the projections as above.
Instead of a transversality condition assume now that the reduced
tangential cohomology $\bar{H}_\tang^* (M)$ of~$M$ is finite
dimensional. Then the Lefschetz number
\begin{displaymath}
  L (\phi^t) := \sum_\kappa (- 1)^\kappa \Tr \bigl( \phi^{t *} : \bar{H}^{\kappa}_\tang (M) \to \bar{H}^{\kappa}_\tang (M) \bigr)
\end{displaymath}
is defined for each~$t$. By the K\"unneth formula $L (\phi^t) $
defines a smooth function~$L (\phi)$ on~$\R_{> 0}$. In this situation,
cf.~Section~\ref{sec:coinform}, one has

\begin{thm} 
  \label{thm:smoothintersection} 
  The intersection product $\cl{\phi} \bullet \cl{\pi}$ is defined.
  Its push forward $\pr_* (\cl{\phi} \bullet \cl{\pi})$ is the
  distribution on~$\R_{> 0}$ which is associated with the function~$L
  (\phi)$.
\end{thm}

To prove this result, besides the K\"unneth theorem, the Poincar\'e
duality isomorphism for~$M$ is the essential tool.

Now suppose in addition that the graphs of $\phi$ and~$\pi$ intersect
foliatedly transversally. Then we can combine the
Theorems~\ref{thm:intersectionproduct}
and~\ref{thm:smoothintersection} to get our version of the dynamical
Lefschetz formula.

\begin{cor} 
  \label{cor:dynlefschetz} 
  For each $t > 0$ we have
  \begin{displaymath}
    L (\phi^t) = \sum_a \sgn \det (\id - \tang_{a} \phi^{t_a}) \frac{1}{\bigl\vert \det \bigl( \id - \trans_a \phi^t \bigr) \bigr\vert}
  \end{displaymath}
  and
  \begin{displaymath}
    \sum_{t / l (\gamma) \in \Z} l (\gamma) \sgn \det (\id - \tang_{b_\gamma} \phi^{t}) \frac{1}{\bigl\vert \det \bigl( \id - \overline{\trans}_{b_\gamma} \phi^{t} \bigr) \bigr\vert} = 0 .
  \end{displaymath}
\end{cor}

\subsection{{}} 
\label{sec:introe}
If $M$ is foliated by points a formula of Guillemin and
Sternberg~\cite[Chapter~VI, p.~311]{MR58:24404} may be read as a
dynamical Lefschetz formula. Here the trace of~$\phi^*$ is defined by
a certain pull back process for distributions.

If $M$ is connected, oriented and foliated as one leaf then
$\bar{H}^{*}_\tang (M) = H^* (M)$ is finite dimensional. By the
homotopy property of de~Rham cohomology $L (\phi)$ is constant and
equal to the Euler characteristic of~$M$. Assume transversal
intersection of the graphs in the classical sense. Then the zeroes
of~$X$ are non degenerate and the Hopf index formula can be
interpreted as a dynamical Lefschetz formula.

Apart from these all results on dynamical Lefschetz formulas are
concerned with fixed point free flows.

For arbitrary leaf dimension the right hand side of
Equation~\eqref{eq:intersectionproduct} appears if the formula of
Guillemin and Sternberg is applied to the exterior bundle of the dual
of~$\tang M$. With varying generality, this is written up by
Guillemin~\cite{MR56:6758}, Deninger~\cite{MR1873325} and Deninger,
Singhof~\cite{MR2002k:58054}. The crucial point now is to find
conditions which justify to pass from the alternating sum of traces on
tangential forms to the alternating sum of traces on tangential
cohomology.

Several examples where this passage is possible for foliated manifolds
with leaves of codimension one are given in the literature:
Guillemin~\cite{MR56:6758} considers the Selberg trace formula.
\'Alvarez L\'opez, Kordyukov~\cite{MR1882768}, Deninger,
Singhof~\cite{MR2002k:58054} and Lazarov~\cite{MR2001i:58044} prove
dynamical Lefschetz formulas for foliated manifolds with a bundle-like
metric.

A dynamical Lefschetz formula does not necessarily hold for foliated
manifolds which do not admit a transversal Riemannian metric. A
counterexample is given by Deninger and Singhof
in~\cite{MR2002c:58055}.

The significance of dynamical Lefschetz formulas in a wider context is
discussed by Deninger~\cite{MR99g:11084,MR2001g:11102} and
Juhl~\cite{MR2001j:37046}.

\begin{ackn}
  I thank Prof.~Dr.~C.~Deninger who motivated this work and supported
  its completion.
\end{ackn}

\section{Foliated Submanifolds and Transversality}
\label{sec:subtrans}

Let $M$ be a foliated manifold. For convenience we write $M^{p, q}$ to
indicate that the leaves of~$M$ have dimension~$p$ and the local
transversal manifolds have dimension~$q$. The pair~$(p, q)$ is called
the foliated dimension of~$M$.

\begin{defn} 
  \label{defn:subman} 
  A subset $S \subseteq M^{p, q}$ is a \emph{foliated submanifold}
  of~$M$ if for every $a \in S$ there is a foliated chart
  \begin{displaymath}
    \phi : U \to U' \times V'
  \end{displaymath}
  around~$a$ with $U' \subseteq \R^r \times \R^{p - r}$ and $V'
  \subseteq \R^s \times \R^{q - s}$, such that
  \begin{displaymath}
    \phi (S \cap U) = \bigl( U' \cap (\R^r \times \{0\}) \bigr) \times \bigl(V' \cap (\R^s \times \{0\} ) \bigr) .
  \end{displaymath}
\end{defn}

The notion foliated codimension will be used with the obvious
meaning.

\begin{exmp} 
  \label{exmp:graph} 
  Let $ f : M \to N $ be a foliated map. Then the graph~$\Gamma_f$
  of~$f$ is a closed foliated submanifold of $M \times N$, which is
  isomorphic to~$M$.
\end{exmp}

\begin{defn} 
  \label{defn:regtrans} 
  Let $ f : M \to N $ be a foliated map and let $S \subseteq N$ be a
  foliated submanifold. The map~$f$ is \emph{foliatedly transversal
    over~S} if
  \begin{align*} 
    \tans_a f (\tans_a M) + \tans_{f (a)} S & = \tans_{f (a)} N, &
    \tang_a f (\tang_a M) + \tang_{f (a)} S & = \tang_{f (a)} N
  \end{align*}  
  for all $a \in f^{-1} (S)$.
\end{defn}

Carrying out the obvious modifications in the proof for the smooth
case we have the following result concerning transversality.

\begin{prop} 
  \label{prop:transsub} 
  Let the foliated map $ f : M \to N $ be foliatedly trans\-versal
  over the foliated submanifold $S \subseteq N$ of foliated
  codimension $(t, u)$. If $f^{-1} (S) \ne \emptyset$, then
  \begin{displaymath}
    f^{-1} (S) \subseteq M
  \end{displaymath}
  is a foliated submanifold of foliated codimension $(t, u)$. \qed
\end{prop}

\section{Foliated Riemannian Connections}
\label{sec:riemconn}

A foliated fiber bundle is a fiber bundle in the category of
foliated manifolds. We call a foliated fiber bundle transversal
if the typical fiber is foliated by points.

Let $M$ be a foliated manifold with foliated tangent bundle
\begin{displaymath}
  \pi : \tans M \to M.
\end{displaymath} 
Suppose that $M$ is equipped with a bundle-like metric denoted by
$\langle - , - \rangle$ and identify $\tans M$ with the Whitney sum $
\tang M \perp \trans M$. Then a smooth vector field~$X$ decomposes
uniquely as
\begin{displaymath}
  X = X_\tang + X_\trans ,
\end{displaymath}
where $X_\tang$, $X_\trans$~resp., is a smooth section of~$\tang M$,
$\trans M$~resp. Moreover, if $X$ is foliated the same is true
for~$X_\trans$.

Let $a \in M$ be a point and let $Y$ be a foliated vector field on a
neighborhood~$U$ of~$a$. For smooth vector fields $X$ and~$Z$ on~$U$
write
\begin{equation} 
  \label{eq:covder} 
  \begin{split} 
    \langle \nabla_X Y , Z \rangle := & \frac{1}{2} \bigl( X \langle Y_\tang , Z_\tang \rangle + Y_\tang \langle Z_\tang , X \rangle - Z_\tang \langle Y_\tang , X  \rangle \\
    & + \langle X , [ Z_\tang , Y_\tang ] \rangle + \langle Y_\tang , [Z_\tang ,  X ] \rangle - \langle Z_\tang , [Y_\tang , X ] \rangle \\
    & + X_\trans \langle Y_\trans, Z_\trans \rangle + Y_\trans  \langle Z_\trans , X_\trans \rangle - Z_\trans \langle Y_\trans , X_\trans \rangle \\
    & + \langle X_\trans , [ Z_\trans , Y_\trans ] \rangle + \langle
    Y_\trans , [Z_\trans , X_\trans ] \rangle - \langle Z_\trans ,
    [Y_\trans ,X_\trans ] \rangle \bigr) .
  \end{split} 
\end{equation} 
Then the assignment
\begin{math} 
  (X, Z) \mapsto \langle \nabla_X Y , Z \rangle
\end{math}
defines a tensor field of type $(0, 2)$ on~$U$. This fact allows the
following

\begin{defn} 
  \label{defn:covder}
  Let $X$ be a smooth vector field on~$U$. Then the \emph{foliated
    covariant derivative of~$Y$ in the direction of~$X$} is the unique
  vector field~$\nabla_X Y$ satisfying equation~\eqref{eq:covder} for
  all vector fields $Z$ on~$U$. If $\xi$ is a tangent vector at~$a$
  then the \emph{foliated covariant derivative~$\nabla_\xi Y$ of~$Y$
    at~$a$ in the direction of~$\xi$} is the value at~$a$ of $\nabla_X
  Y$, where $X$ is any vector field with $X_a = \xi$.
\end{defn}

For a vector space~$V$ and $v \in V$ let
\begin{displaymath}
  \tau_v : \tans_v V \to V
\end{displaymath}
be the canonical identification. If now $\xi$ is a tangent vector
at~$a$ then the vector
\begin{displaymath}
  C (\xi, Y) := \tans_a Y ( \xi ) - \tau_{Y_a}^{-1} ( \nabla_{\xi} Y )
\end{displaymath}
in $\tans_{Y_a} \tans M$ depends only on the value of~$Y$ at~$a$. More
precisely we get

\begin{prop} 
  \label{prop:riemconn} 
  The map
  \begin{align*} 
    C : \tans M \oplus \tans M & \to \tans \tans M \\
    (\xi, z) & \mapsto C (\xi, Y)
  \end{align*} 
  where $Y$ is any local foliated vector field taking the value~$z$,
  is a foliated linear connection on~$M$. \qed
\end{prop}

\begin{defn} 
  \label{defn:riemconn} 
  The foliated linear connection of Proposition~\ref{prop:riemconn} is
  the \emph{foliated Riemannian connection} associated to the
  bundle-like metric on~$M$.
\end{defn}

The main property of the foliated Riemannian connection is that its
parallel transport preserves the metric.

\begin{prop} 
  \label{prop:riemconncomp} 
  Let $ \gamma : I \to M $ be a path. If $\beta_1$ and~$\beta_2$ are
  parallel transports along~$\gamma$ relative to the foliated
  Riemannian connection then $ \frac{d}{d t} \langle \beta_1 (t) ,
  \beta_2 (t) \rangle = 0 $. \qed
\end{prop}

Corresponding to the foliated Riemannian connection on~$M$ there is
the exponential map of its geodesic spray. This is a foliated
map
\begin{displaymath}
  \exp : \mathscr{O} \to M ,
\end{displaymath}
which is defined on an open neighborhood~$\mathscr{O}$ of the zero
section in~$\tans M$.

\section{Foliated Homomorphisms}
\label{sec:folhomo}

Let $E\to{}M$ and $F\to{}N$ be transversal vector bundles and denote
by~$F^*$ the vector bundle dual to~$F$. A foliated homomorphism
$(f,\rho{}):E\to{}F$ consists of a foliated map $f:M\to{}N$ and a
global foliated section~$\rho{}$ of $f^{*}F^{*}\otimes{}E$.  A
foliated homomorphism $(f,\rho{})$ is called proper if $f$ is proper.

\begin{exmp} 
  \label{exmp:exp} 
  Let $M$ be a foliated manifold with a bundle-like metric and let $E
  \to M$ be a transversal vector bundle with a foliated linear
  connection
  \begin{displaymath}
    C : \tans M \oplus E \to \tans E .
  \end{displaymath}
  Let $\xi \in \mathscr{O}$ be a point in the domain of the
  exponential map. Then for each vector $z \in E_{\exp (\xi)}$ there
  exists a unique parallel transport~$\beta_{\xi, z}$ along~$\exp (t
  \cdot \xi)$ with respect to~$C$ satisfying $\beta_{\xi,z}(1)=z$.
  Examining the homogeneous linear differential equation for the
  parallel transport in local coordinates we reveal that the maps
  \begin{align*} 
    P_\xi : E_{\exp (\xi)} & \to E_{\pi (\xi)} \\
    z & \mapsto \beta_{\xi, z} (0)
  \end{align*}  
  fit together to define a foliated section of $\exp^* E^* \otimes
  \pi^* E |_\mathscr{O}$. Hence we get a foliated homomorphism
  \begin{displaymath}
    (\exp,P):\pi^*E|_\mathscr{O}\to{}E. 
  \end{displaymath} 
\end{exmp}

\begin{exmp} 
  \label{exmp:flowhomo} 
  Let $ \pr : \R \times M \to M $ be the projection, where the real
  line is foliated by points. A global foliated flow on~$E$ is
  a foliated homomorphism
  \begin{displaymath}
    (\phi, \rho) : \pr^* E \to E , 
  \end{displaymath} 
  where $\phi$ is a foliated flow on~$M$ and, with the obvious
  notation, $\rho$ satisfies
  \begin{align*} 
    \rho_{a}^0 & = \id, & \rho^{s + t}_a & = \rho^s_{\phi^t (a)} \circ
    \rho^t_{a}
  \end{align*} 
  for all $s, t \in \R$ and $a \in M$.
\end{exmp}

\section{Tangential (Co)homology}
\label{sec:derham}

Let $M^{p, q}$ be a foliated manifold and let $E \to M$ be a
transversal vector bundle of rank~$s$. The tangential (differential)
$k$-forms~$\Omega^k_\tang (M, E)$ on~$M$ with values in~$E$ are the
smooth sections of $\bigwedge^k \tang M^* \otimes E$ equipped with the
$\gCi$-topology. The tangential de~Rham complex of~$M$ with values
in~$E$ consists of the direct sum
\begin{displaymath}
  \Omega^*_\tang (M, E) := \bigoplus_{k \geq 0} \Omega^k_\tang{}(M, E)
\end{displaymath}
together with the tangential exterior differentiation~$d_\tang$, which
is defined locally as follows. Let $x^1, \dots , y^q$ be foliated
coordinates, let $d_\tang x^i$ be the tangential $1$-form which is
dual to~$\frac{\partial}{\partial x^i}$ and let $\sigma_1, \dots ,
\sigma_s$ be a foliated frame for~$E$. Then
\begin{displaymath}
  d_\tang f := \sum_{i = 1}^p \frac{\partial f}{\partial x^i} \: d_\tang x^i 
\end{displaymath} 
for a smooth function~$f$ and
\begin{displaymath} 
  d_\tang \omega := \sum_{I, j} d_\tang a_I^j \wedge d_\tang x^I \otimes \sigma_j
\end{displaymath} 
for a tangential form $\omega = \sum_{I, j} a_I^j \: d_\tang x^I
\otimes \sigma_j$. 

Let $F\to{}M$ be a second transversal vector bundle.
If $\omega\in\Omega^{k}_\tang(M,E)$ and $\tau\in\Omega^{l}_\tang(M,F)$
are given then there is the equality
\begin{displaymath}
  d_\tang(\omega\wedge\tau)=(d_\tang\omega)\wedge\tau+(-1)^{k}\omega\wedge(d_\tang\tau)
\end{displaymath} 
of tangential forms with values in $E\otimes{}F$.

For each compact subset $K \subseteq M$ we denote by
$\Omega^{k}_{\tang{}}(M,E;K)$ the closed vector subspace consisting of
tangential $k$-forms with support contained in~$K$. The space
$\Omega^{k}_{\tang{},c}(M,E)$ of all tangential $k$-forms with compact
support carries its $LF$-topology. Define the pull back of a
tangential form in the obvious way to get contravariant functors
$\Omega^{*}_{\tang{}}$ and~$\Omega^{*}_{\tang{},c}$ from trans\-versal
vector bundles with (proper) foliated homomorphisms to nuclear
graded differential $\R{}$-vector spaces.  The cohomology functors
$H^{*}_{\tang{}}$ and~$H_{\tang{},c}^{*}$ from transversal vector
bundles to topological graded $\R$-vector spaces corresponding to
tangential differential forms are the tangential cohomology and the
tangential cohomology with compact support.

Let the real line~$\R$ be foliated as one leaf and let
$\pr_M:M\times\R\to{}M$ be the projection. Then we have the Poincar\'e
lemma for tangential cohomology, cf.~\cite{MR85a:57016}, i.e.
$\pr_M^*$ defines an isomorphism
\begin{displaymath}
  H^{*}_\tang(M,E)\cong{}H^{*}_\tang(M\times\R,\pr_M^*E).
\end{displaymath}
This implies that tangential cohomology is invariant under certain
homotopies which we call tangential.

Let $\pi:V\to{}M$ be a foliated vector bundle of rank~$r$ the typical
fiber of which is foliated as one leaf. A tangential form
$\omega\in\Omega^*_\tang(V)$ has compact support in the vertical
direction if for each compact $K\subseteq{}M$ the set
$\Supp\omega\cap\pi^{-1}(K)$ is compact. The resulting subcomplex is
denoted by $\Omega^{*}_{\tang,cv}(V)$, the cohomology
$H^{*}_{\tang,cv}(V)$ of this complex is the tangential cohomology
of~$V$ with compact support in the vertical direction. Assume that
$V$ is oriented. Then integration along the fiber~$\pi_*$ induces the
Thom isomorphism, cf.~\cite{MR89h:58184},
\begin{displaymath}
  H^{*}_{\tang,cv}(V)\cong{}H_\tang^{*-r}(M).
\end{displaymath}
The cohomology class $\Phi_\tang(V):=\pi_*^{-1}(1)$ in
$H_{\tang,cv}^{r}(V)$ is the tangential Thom class of~$V$.  Now
fix a Riemannian metric in~$V$. Then for $\epsilon>0$ there exists a
representative of the tangential Thom class with support contained in
\begin{displaymath}
  V(\epsilon):=\{v\in{}V\mid\Vert{}v\Vert<\epsilon\}.
\end{displaymath}

A tangential $E$-current of dimension~$k$ on~$M$ is a continuous
linear form on $\Omega^{k}_{\tang{},c}(M,E)$.  The vector space of all
tangential $E$-currents of dimension~$k$ on~$M$ will be denoted by
$\Omega^{\tang{}}_{k}(M,E)$ and equipped with the weak topology.
Transposition allows to define the tangential boundary
map~$b_{\tang{}}$ and the push forward by a proper foliated
homomorphism. Hence, the tangential currents become a covariant
functor from trans\-versal vector bundles with proper foliated
homomorphisms to topological graded codifferential $\R$-vector spaces.
The homology functor~$H^{\tang{}}_{*}$ associated with the tangential
currents is the tangential homology. Similarly, the tangential
currents with compact support~$\Omega^{\tang{}}_{c,*}$ are dual
to~$\Omega^{*}_{\tang{}}$ and define tangential homology with compact
support~$H^\tang_{c,*}$. We always think of the trivial bundle
$E=\R{}$ if the~$E$ is omitted in the notation.

\begin{exmp}
  \label{exmp:integral}
  Let $M^{p,q}$ be a foliated manifold with a transversal Riemannian
  metric and corresponding transversal volume~$\vol_{\trans{}}$. Let
  $\tang{}M$ be oriented. Then real valued tangential $p$-forms are
  identified with densities and integration defines the closed
  tangential current $\int_{M}-\:\vol_{\trans{}}$ of dimension~$p$
  on~$M$. Moreover, the assignment
  $\omega\mapsto\int_M\omega\wedge-\:\vol_\trans$ induces continuous
  maps
  \begin{align*}
    H^{p-k}_{\tang,c}(M,E^*)&\to{}H^\tang_{c,k}(M,E),&H^{p-k}_\tang(M,E^*)&\to{}H_{k}^\tang{}(M,E).
  \end{align*}
\end{exmp}

Since tangential exterior differentiation is not a topological
homomorphism in general, we consider the reduced tangential
(co-)homology (with compact support) which is the maximal Hausdorff
quotient of the corresponding (co-)homology space. By continuity of
pull back we get new functors, denoted by
$\bar{H}^{*}_\tang$,~$\bar{H}^{*}_{\tang,c}$, $\bar{H}^\tang_{c,*}$
and~$\bar{H}^\tang_{*}$. Use the Hahn-Banach theorem and the theorem
of bipolars to prove

\begin{prop}
  \label{prop:homtocoh}
   The canonical maps
   \begin{align*} 
     \bar{H}^\tang_k(M,E)&\to\bar{H}^{k}_{\tang,c}(M, E)',&\bar{H}^\tang_{c,k}(M,E)&\to\bar{H}^{k}_\tang(M,E)'
   \end{align*}
   are topological isomorphisms. \qed
\end{prop}

\section{K\"unneth Maps} 
\label{sec:kuenneth} 

Let $E\to{}M^{p,q}$ and $F\to{}N^{r,s}$ be transversal vector bundles
and let $E\boxtimes{}F$ be the external tensor product.  Consider the
bilinear map
\begin{gather*}
  {}-\otimes-:\Omega^{*}_\tang(M,E)\times\Omega^{*}_\tang(N,F)\to\Omega^{*}_\tang(M\times{}N,E\boxtimes{}F),\\
  \omega\otimes\tau:=\pr_M^*\omega\wedge\pr_N^*\tau,
\end{gather*}
where $\pr_M$ and~$\pr_N$ are the projections. A bilinear map
\begin{displaymath}
  -\otimes-:\Omega^{\tang}_k(M,E)\times\Omega^{\tang}_l(N,F)\to\Omega^{\tang}_{k+l}(M\times{}N,E\boxtimes{}F)
\end{displaymath}
is defined by the requirement that
\begin{displaymath}
  \langle{}S\otimes{}T,\omega\otimes\tau\rangle=(-1)^{k(r-l)}\langle{}S,\omega\rangle\langle{}T,\tau\rangle
\end{displaymath}
for all $\omega\in\Omega^k_{\tang,c}(M,E)$ and
$\tau\in\Omega^l_{\tang,c}(N,F)$. Both maps restrict to compact
supports and descend to reduced (co)homology. Moreover, by linear
disjointness with respect to $-\otimes-$ the induced maps
\begin{align*}
  \bar{H}^{*}_{\tang(,c)}(M,E)\otimes\bar{H}^{*}_{\tang(,c)}(N,F)&\to\bar{H}^{*}_{\tang(,c)}(M\times{}N,E\boxtimes{}F),\\
  \bar{H}^{\tang}_{(c,)*}(M,E)\otimes\bar{H}^{\tang}_{(c,)*}(N,F)&\to\bar{H}^{\tang}_{(c,)*}(M\times{}N,E\boxtimes{}F)
\end{align*}
are one-to-one. In special cases we get more precise results for the
nuclear cohomology spaces.

\begin{exmp} 
  \label{exmp:kuentrivtrans} 
  Let $M$ be foliated by points. Similar to Theorems~$40.1$ and~$44.1$
  of~\cite{MR37:726} one proves the topological isomorphisms
  \begin{align*}
    \gSec_{(c)}(M,E)\ctens\bar{H}^{*}_{\tang(,c)}(N,F)&\cong\bar{H}^{*}_{\tang(,c)}(M\times{}N,E\boxtimes{}F).\\
  \end{align*}
\end{exmp}

\begin{exmp} 
  \label{exmp:kuenharm} 
  Let $M$ be a compact foliated manifold with bundle-like metric and
  oriented tangential subbundle. Suppose we are given a transversal
  Riemannian metric in~$E$, i.e.~the corresponding section of
  $E^*\otimes{}E^*$ is foliated. Define the tangential star operator
  \begin{displaymath}
    \star_\tang:\Omega^k_\tang(M,E)\to\Omega^{p-k}_\tang(M,E)
  \end{displaymath}
  by the requirement that
  \begin{math}
    \omega\wedge\star_\tang\tau=\langle\omega,\tau\rangle\:\vol_{\tang}
  \end{math}
  for all $\omega\in\Omega^k_\tang(M,E)$ and let
  \begin{displaymath}
    \delta_\tang:=(-1)^{p\cdot*+1}\star_\tang{}d_\tang\star_\tang:\Omega^{*+1}_\tang(M,E)\to\Omega^{*}_\tang(M,E).
  \end{displaymath}
  The tangential Laplacian of~$E$ is the operator
  $\Delta_\tang:=d_\tang\circ\delta_\tang+\delta_\tang\circ{}d_\tang$
  and we write
  \begin{displaymath}
    \mathscr{H}^{k}_\tang(M,E):=\ker\Delta_\tang\cap\Omega^{k}_\tang(M,E) 
  \end{displaymath}
  for the tangentially harmonic $k$-forms. By adjointness with respect
  to the $L^2$-inner product $ \int_M - \wedge \star_\tang - \:
  \vol_\trans$ we find $\mathscr{H}_\tang (M, E) = \ker d_\tang
  \cap \ker \delta_\tang$ hence a map
  \begin{displaymath}
    \mathscr{H}^{*}_\tang(M, E)\to\bar{H}^{*}_\tang(M,E), 
  \end{displaymath}
  which is an inclusion since $\star_\tang \circ \Delta_\tang =
  \Delta_\tang \circ \star_\tang$. It is a consequence of the
  tangential Hodge decomposition
  \begin{displaymath}
    \Omega^*_\tang(M,E)=\mathscr{H}^*_\tang(M,E)\oplus\overline{\im{}d_\tang}\oplus\overline{\im\delta_\tang},
  \end{displaymath}
  cf.~\cite[Corollary~$1.3$]{MR1815391}, that this is a topological
  isomorphism. Assume that $N$ is a compact foliated manifold with
  bundle-like metric and oriented tangential subbundle, too. Let
  $F\to{}N$ be equipped with a transversal Riemannian metric.
  Collecting the signs correctly we find the commutativity of the
  diagram
  \begin{equation} 
    \label{eq:comdia}
    \begin{CD}\Omega^{*}_\tang(M,E)\otimes\Omega^{*}_\tang(N,F)@>{\Delta_\tang^M\otimes\id+\id\otimes\Delta_\tang^N}>>\Omega^{*}_\tang(M,E)\otimes\Omega^{*}_\tang(N,F)\\@VVV@VVV\\\Omega^{*}_\tang(M\times{}N,E\boxtimes{}F)@>{\Delta_\tang^{M\times{}N}}>>\Omega^{*}_\tang(M\times{}N,E\boxtimes{}F)\end{CD}. 
  \end{equation}
  Now consider the Hilbert spaces obtained by completion of
  differential forms with respect to the $L^2$-product, which are
  denoted by a prefixed~$L^2$.  By \cite[Theorem~$2.2$]{MR51:6119},
  the tangential Laplacians define essentially selfadjoint unbounded
  operators, so let $E^M(t)$, $E^N(t)$~resp., be the spectral families
  of $\Delta_\tang^M$, $\Delta_\tang^N$~resp. The operator
  $\Delta_\tang^M\otimes\id+\id\otimes\Delta_\tang^N$ is essentially
  selfadjoint in the complete tensor product
  $L^2\Omega^{*}_\tang(M,E)\ctens{}L^2 \Omega^{*}_\tang(N,F)$ of
  Hilbert spaces. Its spectral family is given by, cf.~\cite[proof of
  Theorem~$8.34$]{MR81e:47001},
  \begin{displaymath}
    E(t)=G(\{z\in\C\mid\Re(z)+\Im(z)\leq{}t\}),
  \end{displaymath}
  where $G$ is the complex spectral family defined by
  \begin{displaymath}
    G(t+is):=E^M(t)\ctens\id\circ\id\ctens{}E^N(s).
  \end{displaymath}
  The isomorphism
  \begin{displaymath}
    L^2\Omega^{*}_\tang(M,E)\ctens{}L^2\Omega^{*}_\tang(N,F)\cong{}L^2\Omega^{*}_\tang(M\times{}N,E\boxtimes{}F) 
  \end{displaymath}
  together with the diagram~\eqref{eq:comdia} allows to identify~$E$
  as the spectral family of~$\Delta_\tang^{M \times N}$. Since the
  tangential Laplacians are nonnegative we get a commutative diagram
  \begin{displaymath}
    \begin{CD}L^2\Omega^{*}_\tang(M,E)\ctens{}L^2\Omega^{*}_\tang(N,F)@>{E^M (0) \ctens E^N (0)}>>L^2\Omega^{*}_\tang(M,E)\ctens{}L^2\Omega^{*}_\tang(N,F)\\@VVV@VVV\\L^2\Omega^{*}_\tang(M\times{}N,E\boxtimes{}F)@>{E(0)}>>L^2\Omega^{*}_\tang(M\times{}N,E\boxtimes{}F)\end{CD}.
  \end{displaymath}
  By the decomposition of tangential forms there is an induced diagram
  \begin{displaymath}
    \begin{CD}\Omega^{*}_\tang(M,E)\ctens\Omega^{*}_\tang(N,F)@>{\Pi^M\ctens\Pi^N}>>\mathscr{H}^{*}_\tang(M,E)\ctens\mathscr{H}^{*}_\tang(N,F)\\@VVV@VVV\\\Omega^{*}_\tang(M\times{}N,E\boxtimes{}F)@>{\Pi}>>\mathscr{H}^{*}_\tang(M\times{}N,E\boxtimes{}F)\end{CD},
  \end{displaymath}
  and we can
  replace~$\mathscr{H}$ by~$\bar{H}$ to find an isomorphism
  \begin{displaymath}
    \bar{H}^{*}_\tang(M,E)\ctens\bar{H}^{*}_\tang(N,F)\cong\bar{H}^{*}_\tang(M\times{}N,E\boxtimes{}F).
  \end{displaymath} 
\end{exmp}

\section{Regularization} 
\label{sec:regular} 

Let $M^{p, q}$ be a foliated manifold with a bundle-like metric in the
foliated tangent bundle $\pi:\tans{}M\to{}M$. Denote by
$\mathcal{V}\tans{}M$ the vertical subbundle of $\tans\tans{}M$ and
let $\upsilon_\trans$ be the transversal density in $\mathcal{V}\tans
M$ induced by the identification $\tans_{\pi(\xi)}
M\cong\mathcal{V}_\xi\tans M$. Let $\tang M$ be oriented and let
$\phi_{\tang, \nu}$ be a representative of the tangential Thom class
of~$M$ with support contained in $\tang{}M(1/\nu)$. Let
$\pr_1:\tans{}M\to\tang{}M$ be the projection and define closed
$p$-forms on $\tans{}M$ by
\begin{displaymath}
  \phi_\nu:=\frac{\pr_1^*\phi_{\tang,\nu}\wedge\nu^q\varrho(\nu\cdot\Vert-\Vert_\trans)}{\pi^*\pi_*\bigl(\pr_1^*\phi_{\tang,\nu}\wedge\varrho(\Vert-\Vert_\trans)\:\upsilon_\trans\bigr)},
\end{displaymath}
where $\varrho\in\gCi(\R)$ is any nonnegative function with
$\varrho(s)=1$ for $s<1/3$, $\varrho(s)=0$ for $s>2/3$. Consider
$M \times N$, $N := M$, and let $\pr_M$ and~$\pr_N$ be the
projections.  Finally, let $\mathscr{U}$ be an open neighborhood of
the zero section in~$\tans{}M$, such that
$(\pi,\exp)\vert_{\mathscr{U}}$ is a foliated embedding.

Suppose now we are given a transversal vector bundle $E\to{}M$ with a
foliated linear connection and recall the foliated homomorphism
$(\exp,P)$ constructed in Example~\ref{exmp:exp}. Let $K\subseteq{}M$
be compact and let $U\supseteq{}K$ be a relatively compact open subset
of~$M$ with closure~$L$. If $\omega\in\Omega_{\tang}^{k}(M,E;K)$ is a
tangential form, then Proposition~\ref{prop:riemconncomp} implies that
for large values of~$\nu$ the support of
$\phi_\nu\wedge(\exp,P)^*\omega$ is contained in a compact subset
of~$\pi^{-1}(L)\cap\mathscr{U}$. Consider this as a form on $M\times{}N$
to define
\begin{displaymath} 
  R'_\nu\omega:=\int_N\phi_\nu\wedge(\exp,P)^*\omega\:\vol_\trans.
\end{displaymath}
By compatibility with the differentials, there are induced linear maps
\begin{displaymath}
  R'_\nu:\bar{H}^k_{\tang}(M,E;K)\to\bar{H}^k_{\tang}(M,E;L).
\end{displaymath}

\begin{lem} 
  \label{lem:formapprox} 
  For each $\omega\in\bar{H}_{\tang}^{k}(M,E;K)$ the sequence
  $R'_\nu\omega$ converges to $(-1)^{pk}\omega$ in
  $\bar{H}^k_{\tang}(M,E;L)$.
\end{lem}

\begin{proof}
  Pull back to the tangent bundle and combine a tangential homotopy
  argument with a familiar estimate.
\end{proof}

Let $S$ be a tangential $E$-current of dimension~$k$ on~$M$ then for
large~$\nu$ Proposition~\ref{prop:riemconncomp} allows to define
\begin{displaymath}
  R_\nu{}S:=\langle{}S_M,\phi_\nu\wedge{}P\rangle
\end{displaymath}
as an element of $\Omega_\tang^{p-k}(U,E^*)$. Again, we get induced
maps
\begin{displaymath}
  R_\nu:\bar{H}^\tang_{k}(M,E)\to\bar{H}_{\tang}^{p-k}(U,E^*).
\end{displaymath}

\begin{lem} 
  \label{lem:approxform} 
  For each $S\in\bar{H}^\tang_k(M,E)$ the sequence
  $\int_{M}R_\nu{}S\wedge-\:\vol_\trans$ converges to $(-1)^{pk}S$ in
  $\bar{H}_\tang^{k}(M,E;K)'$.
\end{lem}

\begin{proof} 
  For $\omega\in\bar{H}_\tang^{k}(M,E;K)$ we find
  \begin{displaymath}
    \int_MR_\nu{}S\wedge\omega\:\vol_\trans=\langle{}S,R'_\nu\omega\rangle
  \end{displaymath}
  and the result follows from Lemma~\ref{lem:formapprox}.
\end{proof}

Conversely, let $S$ be a tangential $E$-current of dimension~$k$ with
support contained in~$K$. Then for large~$\nu$ the form $R_\nu S$ has
support in~$L$ by Proposition~\ref{prop:riemconncomp} and consequently
we get maps
\begin{displaymath} 
  R_\nu:\bar{H}^\tang_{k}(M,E;K)\to\bar{H}^{p-k}_{\tang}(M,E^*;L). 
\end{displaymath}
Similar to the above we have

\begin{lem} 
  \label{lem:approxformc} 
  For $S\in\bar{H}^\tang_{k}(M,E;K)$ the sequence
  $\int_{M}R_\nu{}S\wedge-\:\vol_\trans$ converges to $(-1)^{pk}S$ in
  $\bar{H}_\tang^{k}(M,E)'$. \qed
\end{lem}

This lemma together with Proposition~\ref{prop:homtocoh} yields

\begin{thm} 
  \label{cor:regular} 
  The map 
  \begin{math} 
    \bar{H}^{p-k}_{\tang,c}(M,E^*)\to\bar{H}_\tang^{k}(M,E)'
  \end{math} 
  associated with the integration current has dense image. \qed
\end{thm}

Assume the existence of $\epsilon > 0$ such that
$\tans{}M(\epsilon)\subseteq\mathscr{U}$. Then for large~$\nu$ the
above definition yields maps
\begin{displaymath}
  R_\nu:\bar{H}^\tang_k(M,E)\to\bar{H}_\tang^{p-k}(M,E^*)
\end{displaymath} 
and one obtains Poincar\'e duality.

\begin{thm} 
  \label{thm:regularb} 
  For each~$k$ the bilinear map 
  \begin{displaymath} 
    \int_M-\wedge-\:\vol_\trans:\bar{H}^{p-k}_\tang(M,E^*)\times\bar{H}^{k}_{\tang,c}(M,E)\to\R
  \end{displaymath} 
  is a dual pairing. \qed
\end{thm}

\begin{exmp} 
  \label{exmp:duality} 
  If $M$ is compact an $\epsilon > 0$ as above always exists. In the
  special case where $\bar{H}^{k}_\tang(M, E)$ is finite dimensional
  we get an isomorphism
  \begin{displaymath}
    \bar{H}^{p-k}_\tang(M,E^*)\cong\bar{H}_\tang^{k}(M,E)'.
  \end{displaymath}
\end{exmp}

\section{Product of Homology Classes}
\label{sec:homprod}

Let $M^{p,q}$ be a foliated manifold with a bundle-like metric, let
$\tang{}M$ be oriented and re-use the notations of
Section~\ref{sec:regular}. Let $E$ and~$F$ be transversal vector
bundles over~$M$, where $F$ is equipped with a foliated linear
connection. Let $S\in\bar{H}^\tang_{k}(M,E)$ and
$T\in\bar{H}^\tang_{l}(M,F)$, $k+l\ge{}p$, be given then for
$K\subseteq{}M$ compact 
\begin{displaymath}
  \langle{}S,R_\nu{}T\wedge-\rangle\in\bar{H}_\tang^{k+l-p}(M,E\otimes{}F;K)'
\end{displaymath}
is defined for large values of~$\nu$. If for each cohomology class
$\eta\in\bar{H}_{\tang,c}^{k+l-p}(M,E\otimes{}F)$ the sequence
$\langle{}S,R_\nu{}T\wedge\eta\rangle$ converges we make the following

\begin{defn} 
  \label{defn:homprod} 
  The \emph{intersection product} $S\bullet{}T$ of $S$ and~$T$ is the
  continuous linear form on $\bar{H}_{\tang,c}^{k+l-p}(M,E\otimes{}F)$
  given by
  \begin{displaymath}
    S\bullet{}T:=\lim_{\nu\to\infty}(-1)^{pl}\langle{}S,R_\nu{}T\wedge-\rangle.
  \end{displaymath}
\end{defn}

\begin{exmp} 
  \label{exmp:homprod} 
  Let $\omega\in\bar{H}_\tang^{p-k}(M,E^*)$ and
  $\tau\in\bar{H}^{p-l}_{\tang}(M,F^*)$ be given, then the intersection
  product of the corresponding homology classes is defined and
  \begin{displaymath}
    \Bigl(\int_M\omega\wedge-\:\vol_\trans\Bigr)\bullet\Bigl(\int_M\tau\wedge-\:\vol_\trans\Bigr)=\int_M(\omega\wedge\tau)\wedge-\:\vol_\trans. 
  \end{displaymath}
\end{exmp}

\begin{exmp}
  \label{exmp:intersection}
  Let the closed foliated submanifolds $S$ and~$T$ of~$M$ intersect
  foliatedly transversally, i.e.~the inclusion of~$S$ is foliatedly
  transversal over~$T$, and assume $S\cap{}T\ne\emptyset$. If the
  leaves of~$S$ have dimension~$k$ and those of~$T$ have
  dimension~$l$, then $S\cap{}T$ is a foliated submanifold with leaves
  of dimension $k+l-p$, cf.~Proposition~\ref{prop:transsub}. The
  transversal Riemannian metric on~$M$ induces a transversal
  Riemannian metric on each submanifold and we get transversal
  Riemannian volumes on $S$, $T$ and~$S\cap{}T$.  Assume $\tang{}S$
  and~$\tang{}T$ are oriented, then there is an intrinsic
  proceeding~\cite[$(24.13.14)$]{MR84a:57021} to deduce an orientation
  for~$\tang(S\cap{}T)$. Consequently the closed foliated submanifolds
  define tangential homology classes, again denoted by $S$, $T$
  and~$S\cap{}T$.  A tangential homotopy argument followed by some
  estimation shows that the intersection product $S\bullet{}T$ is
  defined and
  \begin{displaymath}
    S\bullet{}T=(-1)^{k(p-k)+l(p-l)}h\cdot{}S\cap{}T,
  \end{displaymath}
  where $h\in\bar{H}^0_\tang(S\cap{}T)$ is the closed function
  satisfying
  \begin{displaymath}
    h\cdot\vol_\trans(S\cap{}T)\otimes\vol_\trans(M)=\vol_\trans(S)\otimes\vol_\trans(T).
  \end{displaymath}
  By definition of the~$R_\nu$ we have $S\bullet{}T=0$ if $S$
  and~$T$ have empty intersection.
\end{exmp}

\section{Tangential Coincidence} 
\label{sec:coincid}

Let the foliated manifolds $M^{p, q}$ and~$N^{r, s}$, $p\geq{}r$,
carry bundle-like metrics and let the tangential subbundles be
oriented.  Let $E \to M $ and $F \to N$ be transversal vector bundles
and let
\begin{align*} 
  (f, \rho) : E & \to F , & (g, \sigma) : E^* & \to F^*
\end{align*} 
be foliated homomorphisms such that the graphs $\Gamma_f$
and~$\Gamma_g$ intersect foliatedly transversally.  Give
$\tang(M\times{}N)$ the product orientation and $\trans(M\times{}N)$
the product metric, orient a graph by demanding the graph map to be
orientation preserving and $\Gamma_f\cap\Gamma_g$ as in
Section~\ref{sec:homprod}.  Finally, equip
$\trans(\Gamma_f\cap\Gamma_g)$ with the metric induced by that on
$M\times{}N$ and consider $\rho$, $\sigma$~resp., as global foliated
section of $E\boxtimes{}F^*\vert_{\Gamma_f}$,
$E^*\boxtimes{}F\vert_{\Gamma_g}$~resp. Let
$h_f\in\bar{H}^0_\tang(\Gamma_f)$ and
$h_g\in\bar{H}^0_\tang(\Gamma_g)$ be the closed functions such that
\begin{displaymath}
  \frac{1}{h_f}\:\vol_\trans(\Gamma_f)=\vol_\trans(M)=\frac{1}{h_g}\:\vol_\trans(\Gamma_g)
\end{displaymath}
and let $h \in \bar{H}^0_\tang(\Gamma_f\cap\Gamma_g)$ be defined by
\begin{displaymath}
  h\cdot\vol_\trans(\Gamma_f\cap\Gamma_g)\otimes\vol_\trans(M\times{}N)= \vol_\trans(\Gamma_f)\otimes\vol_\trans(\Gamma_g).
\end{displaymath}

\begin{defn} 
  \label{defn:tancoin} 
  The \emph{tangential coincidence} of $(f,\rho)$ and $(g,\sigma)$ is
  the reduced homology class
  \begin{displaymath}
    \Coin_\tang\bigl((f,\rho),(g,\sigma)\bigr):=\int_{\Gamma_f\cap\Gamma_g}\bigl(\frac{\rho}{h_f}\wedge\frac{\sigma}{h_g}\bigr)\wedge-\:h\cdot\vol_\trans
  \end{displaymath}
  in $\bar{H}^\tang_{p-r}(M\times{}N)$.
\end{defn}

Notice that the definition is arranged, such that by
Example~\ref{exmp:intersection}
\begin{displaymath}
  \Coin_\tang\bigl((f,\rho),(g,\sigma)\bigr)=\bigl(\frac{\rho}{h_f}\wedge\Gamma_f\bigr)\bullet\bigl(\frac{\sigma}{h_g}\wedge\Gamma_g\bigr)
\end{displaymath}
if the bundles $E$ and~$F$ are equipped with foliated linear
connections.

\begin{exmp} 
  \label{exmp:lefschetz} 
  If $E=F$ and $(g,\sigma)=\id$ then $\Gamma_f\cap\Gamma_g$ consists
  of the images of the fixed points of~$f$ in the diagonal and the
  tangential coincidence is the distribution
  \begin{displaymath}
    \sum_{f(a)=a}\sgn\det(\id-\tang_af)\frac{\Tr(\rho_{a})}{\vert\det(\id-\trans_af)\vert}\cdot\delta_{(a,a)},
  \end{displaymath}
  where $\delta_{(a,a)}$ is the Dirac distribution at~$(a,a)$.
\end{exmp}

\begin{exmp} 
  \label{exmp:trace} 
  Let again $E=F$ and $(g,\sigma)=\id$ and assume the foliated
  manifold~$M$ to be compact and foliated by points. Choose
  any linear connection in~$E$ then
  \begin{displaymath}
    \int_M\Tr\Bigl(\Delta^*R_\nu\bigl(\frac{\rho}{h_f}\wedge\Gamma_f\bigr)\Bigr)\;\vol_\trans \xrightarrow{\nu \to \infty} \left\langle \Coin_\tang \bigl( (f, \rho), \id \bigr), 1 \right\rangle
  \end{displaymath}
  with $\Delta := \Gamma_{\id}$ by Example~\ref{exmp:intersection}.
  The left hand side is the trace of the smoothing operator defined
  by~$R_\nu(\frac{\rho}{h_f}\wedge\Gamma_f)$ and
  Lemma~\ref{lem:approxformc} suggests to look at the right hand side
  as the trace of the operator defined
  by~$\frac{\rho}{h_f}\wedge\Gamma_f$. Since
  \begin{displaymath}
    \bigl\langle \frac{\rho}{h_f}\wedge\Gamma_f, \omega \otimes \tau \bigr\rangle = \int_M \omega \wedge (f, \rho)^* \tau \; \vol_\trans
  \end{displaymath}
  for $\omega \in \gSec (M, E^*)$ and $\tau \in
  \gSec (M, E)$ this operator is the pull back by~$(f,\rho)$.
\end{exmp}

\begin{exmp} 
  \label{exmp:flowrestriction} 
  Let the positive real numbers $\R_{> 0}$ be foliated by points and
  equipped with the standard metric, let $M = \R_{> 0} \times
  N$, let
  \begin{displaymath}
    g= \pr_N : \R_{> 0} \times N \to N
  \end{displaymath}
  be the projection, let $E = \pr_N^* F$ and let $\sigma$ be the
  identity in each fiber. If
  \begin{displaymath}
    (f, \rho) = (\phi, \rho)|_M
  \end{displaymath}
  is the restriction of a global foliated flow on~$F$, then the
  submanifold~$\Gamma_f\cap\Gamma_g$ has components
  \begin{displaymath}
    \R_{> 0} \times \{a\} \times \{a\} ,
  \end{displaymath}
  where $a$ is a fixed point of~$\phi$, and
  \begin{displaymath}
    \{\nu l (\gamma)\} \times \{(b,b)\in{}N\times{}N\mid{}b\in\gamma\}, 
  \end{displaymath}
  $\nu = 1, 2 \dots$, where $\gamma$ is a periodic orbit with least
  positive period~$l (\gamma)$.  Let $ \pr_{\R_{> 0}} : M \times N \to
  \R_{> 0} $ be the projection and assume $N$ to be compact then $
  \pr_{\R_{>0}, *} $ maps the tangential coincidence to the
  distribution on~$\R_{> 0}$ given by
  \begin{multline*} 
    \sum_a \sgn \det (\id - \tang_{a} \phi^{t_a}) \int_{\R_{> 0}} \frac{\Tr ( \rho^t_{a} )}{\bigl\vert \det \bigl( \id - \trans_a \phi^t \bigr) \bigr\vert} \cdot - \: d t \\
    + \sum_\gamma l (\gamma) \sum_{\nu = 1}^\infty \sgn \det (\id -
    \tang_{b_\gamma} \phi^{\nu l (\gamma)}) \frac{\Tr \bigl( \rho^{\nu
        l (\gamma)}_{b_\gamma} \bigr)}{\bigl\vert \det \bigl( \id -
      \overline{\trans}_{b_\gamma} \phi^{\nu l (\gamma)} \bigr)
      \bigr\vert} \cdot \delta_{\nu l (\gamma)} ,
  \end{multline*}
  where $ t_a > 0 $ and $b_\gamma \in \gamma$ are arbitrary and
  \begin{displaymath}
    \overline{\trans}_{b_\gamma} \phi^{\nu l (\gamma)} : \trans_{b_\gamma} N / \R \cdot \overline{X}_{\phi, b_\gamma} \to \trans_{b_\gamma} N / \R \cdot \overline{X}_{\phi, b_\gamma}
  \end{displaymath}
  is the endomorphism induced by $\trans_{b_\gamma} \phi^{\nu
    l(\gamma)}$.
\end{exmp}

\section{Coincidence Formula}
\label{sec:coinform}

Let $M_1$ be a foliated manifold which is foliated by points and let
$M_2^{p, q_2}$ and $N^{p, s}$ be compact foliated manifolds with
bundle-like metrics and oriented tangential subbundles. Let $ E_2\to
M_2$ and $F\to N$ be transversal vector bundles with foliated linear
connections and transversal Riemannian metrics, identify the bundles
with their duals, write
\begin{math}
  M := M_1 \times M_2
\end{math}
and $E$ for the bundle $\pr_{M_2}^* E_2$ over~$M$. Let
\begin{math}
  (f, \rho): E \to F
\end{math}
be a foliated homomorphism. Now
assume the reduced cohomology spaces $\bar{H}_\tang^{*} (M_2, E_2)$
and $\bar{H}_\tang^{*} (N, F)$ to be of finite dimension. Then we can
define a foliated homomorphism of trivial transversal vector bundles
\begin{displaymath}
  (f, \rho)^\kappa : M_1 \times \bar{H}_\tang^{\kappa} (M_2, E_2) \to \{\pt\} \times \bar{H}_\tang^{\kappa} (N, F)
\end{displaymath} 
by giving the section
\begin{displaymath}
  (f,\rho)^\kappa_{a} := \left\{
    \begin{aligned}  
      \bar{H}_\tang^{\kappa} (N, F) & \to \bar{H}_\tang^{\kappa} (M_2, E_2) \\
      \tau & \mapsto j_a^*\circ(f, \rho)^* \tau
    \end{aligned}  
  \right\} ,
\end{displaymath} 
where $j_a:M_2\to{}M$ is the inclusion map opposite~$a$. This section is
foliated, i.e.~smooth, by the K\"unneth theorem
\begin{displaymath} 
  \gCi (M_1) \otimes \bar{H}_\tang^{*} (M_2, E_2)\cong\bar{H}_\tang^{*} (M, E).
\end{displaymath}
Now let $(g,\sigma):E\to F$ be a
second foliated homomorphism such that the graphs $\Gamma_f$
and~$\Gamma_g$ intersect foliatedly transversally. With the duality of
Example~\ref{exmp:duality} we consider $(g, \sigma)^{p - \kappa}$ as
a foliated homomorphism
\begin{displaymath}
  (g,\sigma)^{p-\kappa}:M_1 \times \bar{H}_\tang^{\kappa} (M_2, E_2)^* \to \{\pt\} \times \bar{H}_\tang^{\kappa} (N, F)^* .
\end{displaymath}
Having chosen a Riemannian metric on~$M_1$ the tangential coincidence
\begin{equation}
  \label{eq:proofc}  
  \Coin_\tang \bigl( (f, \rho)^\kappa, (g, \sigma)^{p - \kappa} \bigr) = \int_{M_1} \bigl((f,\rho)^\kappa \wedge (g,\sigma)^{p - \kappa}\bigr) \cdot - \: \vol_\trans 
\end{equation} 
is a distribution on $M_1 \times \{\pt\} \cong M_1$, which is
associated with a function.

\begin{thm} 
  \label{thm:coinform} 
  Let $ \pr_{M_1} : M \times N \to M_1 $ be the projection, then
  \begin{displaymath}
    \pr_{M_1,*}\Coin_\tang\bigl((f,\rho),(g, \sigma)\bigr)=\sum_\kappa(-1)^\kappa\Coin_\tang\bigl((f,\rho)^\kappa,(g,\sigma)^{p-\kappa}\bigr).
  \end{displaymath}
  In particular, the left hand side is the distribution associated
  with a function on~$M_1$.
\end{thm}

\begin{proof} 
  For explicit calculation let $(\omega^\kappa_1, \dots ,
  \omega^\kappa_{\alpha_\kappa})$ be a basis of
  $\bar{H}_\tang^{\kappa} (M_2, E_2)$ and let $(\omega^{\kappa, 1},
  \dots , \omega^{\kappa, \alpha_\kappa})$ be the basis of
  $\bar{H}_\tang^{p - \kappa} (M_2, E_2)$ satisfying
  \begin{displaymath}
    \int_{M_2} \omega^{\kappa, i_\kappa} \wedge \omega^\kappa_{j_\kappa} \: \vol_\trans = \delta_{j_\kappa}^{i_\kappa} ,
  \end{displaymath}
  cf.~Example~\ref{exmp:duality}. Similarly, let $(\tau^\kappa_1,
  \dots , \tau^\kappa_{\beta_\kappa})$ and $(\tau^{\kappa, 1}, \dots ,
  \tau^{\kappa, \beta_\kappa})$ be dual bases of
  $\bar{H}_\tang^{\kappa} (N, F)$ and $\bar{H}_\tang^{p - \kappa} (N,
  F)$ with
  \begin{displaymath}
    \int_N \tau^{\kappa, k_\kappa} \wedge \tau^\kappa_{l_\kappa} \: \vol_\trans = \delta_{l_\kappa}^{k_\kappa} .
  \end{displaymath}
  The functions $x^{\kappa, j_\kappa}_{l_\kappa} \in \gCi (M_1)$
  should be the coefficients in the representation
  \begin{displaymath}
    (f, \rho)^* \tau^\kappa_{l_\kappa} = \sum_{j_\kappa} x^{\kappa, j_\kappa}_{l_\kappa} \otimes \omega^\kappa_{j_\kappa}
  \end{displaymath}
  and we let $y^{\kappa, k_\kappa}_{i_\kappa} \in \gCi (M_1)$ be
  defined by
  \begin{displaymath}
    (g, \sigma)^* \tau^{\kappa, k_\kappa} = \sum_{i_\kappa} y^{\kappa, k_\kappa}_{i_\kappa} \otimes \omega^{\kappa, i_\kappa} .
  \end{displaymath}
  Then we get
  \begin{equation}
    \label{eq:proofb} 
    (f,\rho)^\kappa \wedge (g,\sigma)^{p - \kappa}= \sum_{i_\kappa, k_\kappa} x^{\kappa, i_\kappa}_{k_\kappa} y^{\kappa, k_\kappa}_{i_\kappa} . 
  \end{equation} 
  By the K\"unneth isomorphism
  \begin{displaymath}
    \bigoplus_\kappa \gCi_c (M_1) \otimes \bar{H}_\tang^{p - \kappa} (M_2, E_2) \otimes \bar{H}_\tang^{\kappa} (N, F) \to \bar{H}^{p}_{\tang,c} (M \times N, E \boxtimes F)
  \end{displaymath}
  we find the equations
  \begin{displaymath}
    \begin{split} 
      \frac{\rho}{h_f} \wedge \Gamma_f & = \int_{M \times N} \sum_\kappa \sum_{j_\kappa, k_\kappa} (- 1)^{p (p - \kappa)} x^{\kappa, j_\kappa}_{k_\kappa} \otimes \omega^\kappa_{j_\kappa} \otimes \tau^{\kappa, k_\kappa} \wedge - \: \vol_\trans , \\
      \frac{\sigma}{h_g} \wedge \Gamma_g & =
      \int_{M \times N} \sum_\kappa \sum_{i_\kappa, l_\kappa} (-
      1)^\kappa y^{\kappa, l_\kappa}_{i_\kappa} \otimes
      \omega^{\kappa, i_\kappa} \otimes \tau^{\kappa}_{l_\kappa}
      \wedge - \: \vol_\trans
    \end{split}
  \end{displaymath}
  of reduced tangential homology classes, where the functions~$h_f$
  and~$h_g$ are as in Definition~\ref{defn:tancoin}.
  Example~\ref{exmp:homprod} shows that
  \begin{equation}
    \label{eq:proofa} 
    \pr_{M_1, *} \Coin_\tang \bigl( (f, \rho), (g, \sigma) \bigr) = \sum_\kappa (- 1)^\kappa \int_{M_1} \sum_{i_\kappa, k_\kappa} x^{\kappa, i_\kappa}_{k_\kappa} y^{\kappa, k_\kappa}_{i_\kappa} \cdot - \: \vol_\trans . 
  \end{equation} 
  Now combine the equations \eqref{eq:proofa}, \eqref{eq:proofb}
  and~\eqref{eq:proofc} to obtain the result.
\end{proof}

\begin{rem} 
  \label{rem:nottrans} 
  If we drop the assumption, that the graphs intersect foliatedly
  transversally, the equality
  \begin{displaymath} 
    \pr_{M_1, *} \Bigl[ \bigl(\frac{\rho}{h_f} \wedge \Gamma_f \bigr) \bullet \bigl( \frac{\sigma}{h_g} \wedge \Gamma_g\bigr)
    \Bigr] =  \sum_\kappa (- 1)^\kappa \Coin_\tang \bigl( (f,
    \rho)^\kappa, (g, \sigma)^{p - \kappa} \bigr) 
  \end{displaymath} 
  is still valid.
\end{rem}

\begin{exmp}
  \label{exmp:lefschetzformula} 
  If $M_1 = \{\pt\}$, $E_2 = F$, $(f_2, \rho_2) = \id$ and $r=0$ we
  recover from Example~\ref{exmp:lefschetz} the classical Lefschetz
  trace formula
  \begin{displaymath}
    \sum_{f(a) = a} \sgn \det (\id - \tans_a f) \Tr (\rho_{a}) = \sum_\kappa (- 1)^\kappa \Tr \bigl( (f, \rho)^\kappa \bigr).
  \end{displaymath}
\end{exmp}

\begin{exmp} 
  \label{exmp:flowrestrictioncont} 
  Let $(f, \rho)$ be the restriction to $\R_{> 0} \times N$ of a
  global foliated flow and let $(g, \sigma)$ be the projection as in
  Example~\ref{exmp:flowrestriction}. Then with the notation
  introduced above for any $t > 0$ we have
  \begin{multline*} 
    \sum_a \sgn \det (\id - \tang_{a} \phi^{t_a}) \frac{\Tr ( \rho^t_{a} )}{\bigl\vert \det \bigl( \id - \trans_a \phi^t \bigr) \bigr\vert} \\
    = \sum_\kappa (- 1)^\kappa \Tr \bigl( (\phi^t, \rho^t)^* :
    \bar{H}^{\kappa}_\tang (N) \to \bar{H}^{\kappa}_\tang (N) \bigr)
  \end{multline*} 
  and
  \begin{displaymath}
    \sum_{t / l (\gamma) \in \Z} l (\gamma) \sgn \det (\id - \tang_{b_\gamma} \phi^{\nu l (\gamma)}) \frac{\Tr ( \rho^{t}_{b_\gamma} )}{\bigl\vert \det \bigl( \id - \overline{\trans}_{b_\gamma} \phi^{t} \bigr) \bigr\vert} = 0 .
  \end{displaymath}
\end{exmp}

\end{document}